\documentstyle{article}

\renewcommand{\refname}

\font\tenmsb=msbm10 scaled \magstep1
\font\sevenmsb=msbm7 scaled \magstep1
\font\fivemsb=msbm5 scaled \magstep1

\newfam\msbfam
\textfont\msbfam=\tenmsb
\scriptfont\msbfam=\sevenmsb
\scriptscriptfont\msbfam=\fivemsb
\def\Bbb#1{{\fam\msbfam\relax#1}}

\newtheorem{th}{Theorem}
\newtheorem{lem}{Lemma}
\newtheorem{cor}{Corollary}
\begin{document}

\def\R{{\Bbb R}}
\def\C{{\Bbb C}}
\def\N{{\Bbb N}}
\def\Z{{\Bbb Z}}
\def\Im{{\rm Im}\,}
\def\Re{{\rm Re}\,}
\def\e{\epsilon}
\def\a{\alpha}
\def\d{\delta}
\def\l{\lambda}
\def\Subset{\subset\subset}
\def\e{\varepsilon}
\def\supp{{\rm supp}\,}
\def\D{ \noindent $ \Box$}
\def\KK { \par \hspace{10cm} \rule{7pt}{7pt}\par }

\title{}
\author{}
\maketitle
{\rm   УДК 517.518.6}
       \hfill\raisebox{0.0ex}[1.ex][3ex]{}
     {\rm MSC 42A74 (30D35)}
     \newline
      \bf
 Meromorphic almost periodic functions.

\author{ N.\,Parfyonova, S.Ju.\,Favorov.}

\date{}

\maketitle

A function $f\in C(\R)$ is called {\it almost periodic} (a.p.)
if for every $\e >0$ the set of $\e$-almost periods
$$
E_{\e}(f)=\{\tau\in \R:|f(x+\tau)-f(x)|<\e,\quad\forall x\in \R\}
$$
is relatively dense in $\R$. The latter means that for some $L>0$ and
for every $a\in\R$ the segments $[a,a+L]$ have common points with the set
$E_\e(f)$. An analytic function $f(z)$
on a horizontal strip $S$ is said to be {\it analytic a.p.\,function}
if for every $\e >0$ and every substrip $S'\Subset S$ the set of
$\e,\,S'$-almost periods \footnote{For $S=\{z\in\C:\,a<|\Im z|<b\},\
-\infty\le a<b\le\infty,\ S'=\{z\in\C:\,\a<|\Im z|<\beta\}$ the enclose
$S'\Subset S$ means that $a<\a<\beta<b$.}
$$
E_{\e,S'}(f)=\{\tau\in \R:|f(z+\tau)-f(z)|<\e,\quad \forall z\in S'\}
$$
is relatively dense in $\R$.

Theory of analytic a.p.\,functions was constructed by H.\,Bohr, K.\,Bush,
B.\,Jessen; for its detailed presentation, see \cite{JT, L}.
The further development of this theory is closely connected with the names of
M.G.\,Krein, B.Ja.\,Levin, V.P.\,Potapov, H.\,Tornehave, L.I.\,Ronkin
(\cite{KL}-\cite{R4}).

Let us give the following definition:

A meromorphic function $f(z)$ on a strip $S$ (of finite or infinite width)
is called {\it meromorphic a.p.\,function} on this strip
if for each substrip $S'\Subset S$  and each $\e>0$ the set of
$\e,\,S'$-almost periods
$$
E_{\e,S'}(f)=\{\tau\in \R:\rho(f(z+\tau),f(z))<\e,\quad\forall z\in S'\}
$$
is relatively dense in $\R$; by $\rho$ we denote the spherical metric
on $\overline{\C}=\C\cup\{\infty\}$.

Notice that analytic a.p.\,functions are bounded on each
substrip $S'\Subset S$ and the spherical metric is equivalent to
the Euclidean
one on any bounded set; therefore any analytic a.p.\,function is a
meromorphic a.p.\,function, too. It is easy to see that any uniformly
continuous function on $\overline{\C}$ takes
meromorphic a.p.\,functions to  meromorphic a.p.\,functions.
Hence each linear--fractional transformation takes a meromorphic
a.p.\,function (in particular, an analytic a.p.\,function) to a meromorphic
a.p.\,function.

We show that the class of meromorphic a.p.\,functions is not closed
under the operations of addition and multiplication.
Nevertheless the class of meromorphic
a.p.\,functions inherits many properties of the class of analytic
a.p.\,functions: the uniform limit of meromorphic a.p.\,functions is
a meromorphic a.p.\,function; zeros and poles of meromorphic a.p.\,functions
form almost periodic sets; Bochner's criterion of almost periodicity
is valid for meromorphic a.p.\,functions.  Also, we show that every
meromorphic a.p.\,function is a quotient of analytic a.p.\,functions;
this result is based on the theorems about a.p.\,discrete sets from
\cite{RRF, FRR}. Finally, using methods of \cite{F1, F2},
we give a criterion for a pair of a.p.\,sets to be the zero set
and the pole set of some meromorphic a.p.\,function.

At first we prove the following simple statement.

\begin{th}
\label{uni}
A meromorphic a.p.\,function $f(z)$ on a strip $S$ is uniformly continuous
with respect to the metric $\rho$ on any substrip $S'\Subset S.$
\end{th}

\D \ Take $\e>0$. Let $L$ be a positive number such that each segment
$[a,a+L],\ a\in\R,$ contains an $\e,S'$-almost period of $f$.
A meromorphic function is uniformly continuous with respect to the
spherical metric on every compact set, therefore for some $\delta>0$
and any $z,z'\in S'\cap\{z:|\Re z|<~L+1\}$ such that $|z-z'|<\d$
we have $\rho(f(z),f(z'))<\e$.
If $z,z'$ are arbitrary points in $S'$ and $|z-z'|<\d$, then
there exists $\tau\in E_{\e,S'}(f)$ such that $|\Re(z+\tau)|<L+1,
\quad |\Re(z'+\tau)|<L+1$. We get
$$
\rho(f(z),f(z'))\le \rho(f(z),f(z+\tau))+\rho(f(z+\tau),f(z'+\tau))+
\rho(f(z'+\tau),f(z'))<3\e.
$$
\hfill\rule{7pt}{7pt}

\begin{cor}
\label{sep}
If $f$ is a meromorphic a.p.\,function and $N,\ P$ its zero set and
pole set respectively, then we have
$$
{\rm inf}\{|z-z'|:z\in N\cap S',\,z'\in P\cap S'\}\ge\d(S',f)>0
\quad\forall S'\Subset S.                         \eqno(1)
$$
\end{cor}

\D\ It follows from Theorem \ref{uni} that there are no
sequences $z_n,\ z'_n\in S'$ such that $f(z_n)=0,\ f(z'_n)=\infty$ and
$|z_n-z'_n|\to 0$ as $n\to \infty$.
\hfill \rule{7pt}{7pt}

\medskip
We shall say that some property is valid {\it inside $S$} if it is valid
on any substrip $S'\Subset S$. In particular, we shall say that sets
$N$ and $P$ are {\it separated} inside $S$ if they satisfy (1).

Now we can give examples of meromorphic a.p.\,functions such that
their sum or product is not a meromorphic a.p.\,function.

Take $f_1(z)=\sin{\sqrt2\pi z},\ f_2(z)=1/\sin{\pi z}.$
By Kronecker's Theorem (see for example \cite{L}), for any $\d>0$
there exist $t\in\R$ and $m,\,n\in\Z$ such that $|t-n|<\d,\
|t/\sqrt 2-m|<\d$. The points $z'=n/\sqrt 2$ and $z''=m$ belong to
the set of zeros and the set of poles, respectively, for the product
$f_1(z)f_2(z)=\sin{\sqrt2\pi z}/\sin{\pi z}$. Furthermore, we have
$|n/\sqrt 2 - m| \le |n/\sqrt 2-t/\sqrt2|+
|t/\sqrt2-m|\le(1+1/\sqrt2)\d$. Since the choice of $\d$
is arbitrary, zeros and poles of $f_1(z)f_2(z)$
are not separated and $f_1(z)f_2(z)$ is not a meromorphic a.p.\,function.
By the same argument, the sum $1/\sin{\pi z}+1/\sin{(2\sqrt2-1)\pi z}$
is not a meromorphic a.p.\,function, too.

However the class of meromorphic a.p.\,functions is closed with respect
to the uniform convergence:

\begin{th}
\label{conv}
If a sequence of meromorphic a.p.\,functions $f_n(z)$ on $S$
converges uniformly inside $S$, then the limit function $f(z)$
is a meromorphic a.p.\,function on $S$.
\end{th}

\D\ It can easily be checked that the uniform limit of meromorphic functions
with respect to the spherical metric is also a meromorphic function.
Now let n be large enough and $\tau\in E_{\e,S'}(f_n)$. Then we have
$$
\rho(f(z+\tau),f(z))\le\rho(f(z+\tau),f_n(z+\tau))+\rho(f_n(z+\tau),
f_n(z))+\rho(f_n(z),f(z))<3\e.
$$
Hence $f(z)$ is a meromorphic a.p.\,function on $S$.
\hfill\rule{7pt}{7pt}

\medskip
The following result is very useful for the sequel.

\begin{th}
\label{Boch}\ (Bochner's criterion)
The following conditions are equivalent:

\noindent(i) $f(z)$ is a meromorphic a.p.\,function on $S$;

\noindent(ii) for any sequence of real numbers $t_n$ there exists
a subsequence $t_n'$ such that the sequence of functions $f(z+t'_n)$
converges uniformly inside $S$.
\end{th}

\D\ The proof of this theorem is the same as the proof of Bochner's
criterion
for a.p.\,functions on the axis (see, for example, \cite{L}).
\hfill\rule{7pt}{7pt}

\medskip
Note that from Theorem \ref{conv} it follows that the sequence
of functions $f(z+t'_n)$ converges to a meromorphic a.p.\,function inside $S$.

Let us prove the following statement.

\begin{th}
\label{bound}
Let $f(z)$ be a meromorphic a.p.\,function on $S$ and let $P_r$ be the union
of the disks of radius $r$ with the centers at the poles of $f(z)$;
then $f(z)$ is bounded on the set
$S'\setminus P_r$ for any substrip $S'\Subset S$.
\end{th}

\D \ Assume the contrary. Then there exists a sequence of points
$z_n=x_n+iy_n\in S'\setminus P_r$ such that $f(z_n)\to\infty$.
Taking into account Theorem \ref{Boch}, we may assume that the sequence
of functions $f(z+x_n)$ converges uniformly with respect to the
spherical metric to a meromorphic a.p.\,function $g(z)$ inside $S$
and the points $iy_n$ converge to $iy_0\in S$.
\
By uniform continuity of $f(z)$ we get that the functions $f(z+z_n-iy_0)$
also converge uniformly to $g(z)$. In particular, the point $iy_0$ is
a pole of $g(z)$. Suppose $B_0$ is a disk of radius $r'<r$
with the center $iy_0$ such that $\overline{B_0}\subset S$ and
$g(z)$ has no zeros on $\overline{B_0}$ and no poles on $\partial B_0$;
then $|g(z)|\ge\a>0$ for $z\in\overline{B_0}$. Hence the uniform convergence
of the functions $f(z+z_n-iy_0)$ to $g(z)$ with respect to the spherical
metric
implies the convergence of the functions $1/f(z+z_n-iy_0)$ to $1/g(z)$
with respect to the Euclidean metric. Hurwitz' Theorem yields that
the functions $f(z+z_n-iy_0)$ have poles $w_n$ in $B_0$ for $n$
large enough. So the points $w_n+z_n-iy_0$ are poles of $f(z)$ and
$|w_n+z_n-iy_0-z_n|<r$. This contradicts the choice of $z_n$.
\hfill\rule{7pt}{7pt}

\medskip
Let $N_r$ be the union of the disks of radius $r$ with the centers at the
zeros of $f$. Applying Theorem \ref{bound} to the function $1/f$,
we get the inequality
$$
\inf \{|f(z)|:z\in S'\setminus N_r\}>0.
$$

\begin{cor}
\label{anal}
If a meromorphic a.p.\,function has no poles on $S$,
then it is an analytic a.p.\,function on $S$.
\end{cor}

\D\ By the previous theorem, the function $|f(z)|$ is bounded on each
substrip $S'\Subset S$. At the same time the spherical metric is equivalent
to the Euclidean one on any bounded set.
\hfill\rule{7pt}{7pt}

\medskip
Now we can give a simple condition for the product of
meromorphic a.p.\,functions to be a meromorphic a.p.\,function.

\begin{th}
\label{prod}
Let $f_1(z),\ f_2(z)$ be meromorphic a.p.\,functions.
A necessary and sufficient conditions for the product $f_1(z)f_2(z)$ to be
a meromorphic a.p.\,function is that zeros and poles of this product
be separated inside $S$.
\end{th}

\D\  The necessity follows from Corollary \ref{sep}.
Let us prove the sufficiency. Taking into account Theorem \ref{Boch},
we shall show that the uniform convergence of $f_1(z+t_n)$ to $g_1(z)$ and
$f_2(z+t_n)$ to $g_2(z)$ (with respect to the spherical metric) inside $S$
implies the uniform convergence of the functions $f_1(z+t_n)f_2(z+t_n)$
to the function $g_1(z)g_2(z)$.

First we shall show that the distance between zeros and poles
of the function $g_1(z)g_2(z)$ in each substrip $S'\Subset S$
equals the distance between zeros and poles
of the function  $f_1(z)f_2(z)$.

Suppose $g_1(z')g_2(z')=0,\ g_1(z'')g_2(z'')=\infty ,\
z',\ z''\in S'\Subset S$.
Let $C(z'),\ C(z'')$ be the circles of radius $\delta$ with the centers
at the points $z',\,z''$ respectively such that
no zeros and poles of the functions $g_1(z),\,g_2(z)$ lie
either on these circles or inside these circles, except for the
centers.
It can be assumed also that $C(z')\subset S',\ C(z'')\subset S'$.
Using Hurwitz' Theorem for the functions $f_1(z+t_n)$
and $1/f_1(z+t_n)$  we obtain that for $n$ large enough, each function
$f_1(z+t_n)$ has just the same number of zeros and poles
(with the multiplicity) inside the circles $C(z'),\,C(z'')$
as the function $g_1(z)$ has at the points
$z',\,z''$ respectively. The same assertion is true for the functions
$f_2(z+t_n)$ and $g_2(z)$. Hence the number of zeros of the product
$f_1(z+t_n)f_2(z+t_n)$ inside the circle $C(z')$ is greater than
the number of its poles. Conversely, the number of poles of this product
inside the circle $C(z'')$ is greater than the number of its zeros.
It follows that the function $f_1(z+t_n)f_2(z+t_n)$
has at least one zero inside the circle
$C(z')$ and at least one pole inside the circle $C(z'')$.
Therefore the distance between zeros and poles of the product
$f_1(z)f_2(z)$ is not greater than $|z'-z''|+2\d$ with an arbitrary small $\d$.

On the other hand, the uniform convergence of $g_1(z-t_n)$ to $f_1(z)$ and
$g_2(z-t_n)$ to $f_2(z)$ inside $S$ implies that the distance between
zeros and poles of the product $g_1(z)g_2(z)$ is not greater
than $|w'-w''|+2\d$ for an arbitrary zero $w'\in S'$ and a
pole $w''\in S'$ of the product $f_1(z)f_2(z)$.

Furthermore, suppose our theorem is not true. Then there exist
$\gamma>0$ and sequence of points $z_n=x_n+iy_n \in S'\Subset S$
such that
$$
\rho(f_1(z_n+t_n) f_2(z_n+t_n), g_1(z_n) g_2(z_n))\ge \gamma. \eqno(2)
$$

Note that $g_1(z),\,g_2(z)$ are meromorphic a.p.\,functions.
Hence we can assume without loss of generality that
the functions $g_j(z+x_n)$ converge uniformly inside $S$
to meromorphic a.p.\,functions $h_j(z),\ j=1,2,$ respectively and
the points $iy_n$ converge to a point $iy_0\in S$. It is clear that the
functions $f_j(z+t_n+x_n)$ converge uniformly as $n\to\infty$ to the
same functions $h_j(z),\ 1,2$. Since the functions $f_j(z)$ are
uniformly continuious, we see that $f_j(z+z_n+t_n-iy_0)$ also converge
uniformly inside $S$ to the functions $h_j(z),\ j=1,2$.

If the both functions $h_j(z)$ are finite at the point $iy_0$, then the
sequences of functions $f_j(z+z_n+t_n-iy_0)$
and $g_j(z+z_n-iy_0),\ j=1,2,$ are bounded at this point.
Therefore the sequences of numbers $f_1(z_n+t_n)f_2(z_n+t_n)$
and $g_1(z_n)g_2(z_n)$ have the common limit $h_1(iy_0)h_2(iy_0)$.
This contradicts inequality (2).

If the both functions $h_j(z)$ are non vanishing at the point $iy_0$,
then the sequences of functions $1/f_j(z+z_n+t_n-iy_0)$ and
$1/g_j(z+z_n-iy_0),\ j=1,2,$ are bounded at this point.
Therefore the sequences of numbers $1/[f_1(z_n+t_n)f_2(z_n+t_n)]$
and $1/[g_1(z_n)g_2(z_n)]$ have the common limit
$1/[h_1(iy_0)h_2(iy_0)]$. This also contradicts inequality (2).

Now suppose that the function $h_1(z)$ has a zero of multiplicity $k$
at the point $iy_0$ and the function $h_2(z)$ has a pole of multiplicity
$p\le k$ at the same point. Let $C$ be a circle $\{z:\ |z-iy_0|=\d\}$
such that $C\subset S$ and the functions $h_1,\,h_2$
have neither zeros nor poles on the set $\{z:\ 0<|z-iy_0|\le\d\}$.
Hurwitz' Theorem yields that each function $f_1(z+z_n+t_n-iy_0),\
g_1(z+z_n-iy_0)$ has just $k$ zeros inside $C$ for $n>n_1$
and so has $g_1(z+z_n-iy_0)$. By the same reason,
the functions $f_2(z+z_n+t_n-iy_0),\ g_2(z+z_n-iy_0)$ have
just $p$ poles inside $C$ for $n>n_2$.
Since zeros and poles of the functions $f_1(z)f_2(z)$ and $g_1(z)g_2(z)$
are separated inside $S$, we see that
only zeros of the products $f_1(z+x_n+t_n)f_2(z+x_n+t_n),\
g_1(z+x_n)g_2(z+x_n)$ can be in the disk $\{z:\ |z-iy_0|<\d\}$
for $\d$ small enough.

Further, the moduli of the functions $f_1(z+x_n+t_n),\ f_2(z+x_n+t_n),\
g_1(z+x_n),\ g_2(z+x_n)$ are bounded from above and
bounded away from zero uniformly on the circle $C$ for $n$ large enough.
Hence for all $z\in C$ and $n>n_3$,
$$
|f_1(z+z_n+t_n-iy_0) f_2(z+z_n+t_n-iy_0)-
g_1(z+z_n-iy_0)g_2(z+z_n-iy_0)|>\gamma.
$$
By the Maximum Principle, we obtain that the same assertion is valid for
$z=iy_0$. This also contradicts inequality (2).

The same argument for the functions
$1/[f_1(z_n+t_n)f_2(z_n+t_n)]$ and $1/[g_1(z_n)g_2(z_n)]$
leads to the contradiction with (2) in the case $k<p$.
\hfill \rule{7pt}{7pt}

\medskip
\begin{cor}
\label{quot_1}
The quotient of two analytic a.p.\,functions on $S$
is a meromorphic a.p.\,function on $S$ iff
zeros and poles of this quotient are separated inside $S$.
\hfill\rule{7pt}{7pt}
\end{cor}

For studying zeros and poles of meromorphic a.p.\,functions
we use the concept of divisor.

{\it A divisor $d$} in the domain $G\subset\C$ is a set of pairs
$\{(a_n,k_n)\}$ such that the support of the divisor $\supp d_n=\{a_n\}$
is a discrete set in $G$ without limit points in $G$ and
{\it multiplicities} $k_n$ are numbers from $\Z\setminus\{0\}$.
{\it The divisor $d_f$ of a meromorphic
function $f$} on $G$ is the set $\{(a_n,k_n)\}$ such that $a_n$ are
zeros or poles of $f(z)$ and $|k_n|$ are multiplicities of zeros and poles,
with $k_n>0$ for zeros and $k_n<0$ for poles. A divisor
$\{(a_n,k_n)\}$ can be identified
with the charge concentrated on the set of points
$a_n$ with masses $k_n$; the charge for the divisor of a meromorphic
function $f$ is equal to $(1/2\pi)\triangle\log|f|$, where the Laplace
operator is considered in the sense of distributions.
The sum of divisors $d=\{a_n,k_n\}$, $d'=\{a'_n,k'_n\}$ is a divisor
with the support in $\{a'_n\}\cup\{a_n\}$ and the corresponding
multiplicities; in particular, the multiplicity at the point
$a_m=a'_n\in \supp d\cap\supp d'$ is equal to $k_m+k'_n$.
It is clear that $d_{fg}=d_f+d_g$. Further, let $d=\{(a_n,k_n)\}$ be
a divisor; define
$|d|=\{(a_n,|k_n|):\,(a_n,k_n)\in d\}$,\
$d^{+}=\{(a_n,|k_n|):\,(a_n,k_n)\in d,\ k_n>0\}$,\
$d^{-}=\{(a_n,|k_n|):\,(a_n,k_n)\in d,\ k_n<0\}$.
Note that $|d|=d^{+}+d^{-},\ d+d^{-}=d^{+}$.

A divisor $d$ in a strip $S$ is called {\it almost periodic} if
for any smooth (infinite differentiabled) function $\chi(z)$ 
with the support in $S$ the sum $\sum k_n\chi(a_n+t)$,
i.e.\, the convolution of the charge $d$ with the function $\chi$, 
is an a.p.\,function of $t\in \R$ (see \cite{R4}).

A divisor $d=\{(a_n,k_n)\}$ is {\it positive} if $k_n>0$ for all $n$;
this divisor can be identified with the sequence $\{b_m\}$ of points
$\{a_n\}$ where each $a_n$ appears $k_n$ times.
The sequence $b_m$ in $S$ is said to be {\it almost periodic} if
for any $S'\Subset S$ and $\e>0$ the set
\begin{eqnarray*}
E_{\e,S'}=\{t\in \R:\ {\rm there exists a bijection}\quad \a:\N\to \N
\\ {\rm such\ that }\quad
b_m\in S'\& b_{\a(m)}\in S'\Rightarrow |b_m+t-b_{\a(m)}|<\e\}
\end{eqnarray*}
is relatively dense in $\R$ (see \cite{T1} and, in a special case,
\cite{KL}).
This definition is equivalent to the definition of almost periodicity
for the corresponding positive divisor (see  \cite{RRF,FRR}).

It was proved in  \cite{T1} that the divisor $d_f$ of every analytic
a.p.\,function $f$ is almost periodic. We extend this result to
meromorphic f.p.\,functions:

\begin{th}
\label{div}
Suppose $f$ is a meromorphic a.p.\,function on a strip $S$; then
its divisor $d_f$, divisor of zeros $d_f^{+}$, divisor of poles $d_f^{-}$,
and divisor $|d|$ are almost periodic.
\end{th}
\D\  We need the following simple generalization of Lemma 3.1 from \cite{R4}:

\begin{lem}
Let $g,\,f_n, \, n\in\N,$ be meromorphic a.p.\,functions on a domain
$G\subset \C$. If $\rho(f_n(z),\,g(z))\to 0$ uniformly on compact subset of
$G$,
then the functions $\log|f_n|$ considered as the distributions on $G$
converge to $\log|g|$,  and the charges
$d_{f_n}$ considered as the distributions on $G$ converge to the charge $d_g$.
\end{lem}

\D It suffices to check the convergence of the functions $\log|f_n|$
to $\log|g|$ on a neighborhood of each point $z'\in G$.
If $g(z')\neq\infty$, then $g(z)$ is bounded on some neighborhood $U$
of $z'$. Hence the functions $f_n(z)$ converge uniformly on $U$ to $g(z)$
with respect to the Euclidean metric. Using lemma 3.1 from \cite{R4}
we get the convergence of the distributions $\log|f_n|$ to $\log|g|$ on $U$.
If $g(z')=\infty$, then the functions $1/f_n(z)$ converge uniformly
on some neighborhood of $z'$ to $1/g(z)$ with respect to the Euclidean
metric and we can use lemma 3.1 from \cite{R4} again.
Note that $d_{f_n}=(1/2\pi)\triangle\log|f_n|$ and
$d_g=(1/2\pi)\triangle\log|g|$;
since the differentiation keeps the convergence in the sense of
distributions, we obtain the last assertion of the lemma.
\hfill \rule{7pt}{7pt}

\medskip
We continue the proof of the theorem. Let us show that
the convolution
$$
(d_f*\chi)(t)=\sum_nk_n\chi(a_n+t)
$$
of the charge $d_f=\{a_n,\,k_n\}$ with an infinite
differentiated function $\chi(z)$ with the support in $S$
is an a.p.\,function of $t\in \R$.

Let $\{s_n\}$ be a sequence of real numbers.
Taking into account Theorem \ref{Boch},  we may assume that
the functions $f(s_n+z)$ converge, with respect to the spherical metric
uniformly inside $S$, to a meromorphic a.p.\,function $g(z)$.
Let us check that the functions $(d_f*\chi)(s_n+t)$
converge uniformly on $t\in\R$ to the function $(d_g*\chi)(t)$.

Assume the contrary. Then for any $\d>0$ there exists a sequence $t_{n'}\R$
such that
$$
|(d_f*\chi)(s_{n'}+t_{n'})-(d_g*\chi)(t_{n'})|\ge \d. \eqno(3)
$$
As before, it can be assumed that the functions $f(s_{n'}+t_{n'}+z)$ converge,
with respect to the spherical metric uniformly inside $S$, to a
meromorphic a.p.\,function $h(z)$ and so do the functions $g(t_{n'}+z)$.
By the lemma, it follows that the divisors of the functions
$f(s_{n'}+t_{n'}+z)$
converge in sense of distributions to the divisor $d_h$ and so do the
divisors of the functions $g(t_{n'}+z)$. Therefore the functions of $t$\
$(d_f*\chi)(s_{n'}+t_{n'}+t), \ (d_g*\chi)(t_{n'}+t)$ converge to the
same function $(d_h*\chi)(t)$. This contradicts (3).

Now it follows from Bochner's criterion for a.p.\,functions on the axis
(see \cite{L}) that $(d_f*\chi)(t)$ is an a.p.\,function.
Since this statement is true for all smooth
and supported in $S$ functions $\chi(z)$, we see that
$d_f$ is an a.p.\,divisor.

Further, let $\phi(z)$ be a nonnegative smooth function
with the support in a disk $B_r\subset S'\Subset S$ of radius $r$.
Since zeros and poles of $f$ are separated inside $S$, we see that
for $r<r(S')$ and for every $t\in\R$ the support of the function
$\phi(z+t)$ does not contains simultaneously zeros and poles of $f(z)$.
Hence we have $(d^+_f*\phi)(t)=\max\{(d_f*\phi)(t),0\}$. Consequently
the function $(d^+_f*\phi)(t)$ is an a.p.\,function of $t$.
Since any smooth function with support in $S'\Subset S$
is a linear combination of smooth functions
with supports in disks of radius $r<r(S')$, we see that
$d^+_f$ is an a.p.\,divisor. Evidently, $d^-_f=d_f+d^+_f$ and
$|d_f|=d_f+2d^+_f$ are also a.p.\,divisors.
\hfill \rule{7pt}{7pt}

\begin{cor}
\label{div_zer}
Numbers of zeros and poles of meromorphic a.p.\,function
inside the rectangle
$$
\Pi_1(S',t)=\{z\in S': |\Re z-t|<1\},\quad t\in\R,\ S'\Subset S,
$$
are uniformly bounded from above by a constant depending on
$S'$ only.
\end{cor}

\D \ We shall prove that the numbers $|d|(\Pi_1(S',t))$
are bounded uniformly with respect to $t\in\R$. Let $\phi(z)$ be a
nonnegative 
such that $\chi(z)=1$ on the set $\Pi_1(S',0))$. Since the convolution
$(\chi*|d|)(t)$ is an a.p.\,function on $\R$, we see that this convolution
is bounded. Therefore the numbers $|d|(\Pi_1(S',t))\le(\chi*|d|)(t)$
are bounded uniformly with respect to $t\in\R$.
\hfill\rule{7pt}{7pt}

\medskip
The following theorem with Corollary \ref{quot_1}
gives the complete description of meromorphic a.p.\,functions.

\begin{th}
\label{quot_2}
Any meromorphic a.p.\,function $f(z)$ on a strip $S$ is a quotient of two
analytic a.p.\,functions on $S$.
\end{th}

\D\  By Theorem \ref{div}, poles of $f$ form the a.p.\,divisor
$d^-_f$. It follows from \cite{RRF, FRR} that there exists an a.p.\,divisor
$d'$ in $S$ such that $d^{-}+d'$ is the divisor of some a.p.\,analytic
function $h(z)$ on $S$. Then the function $g(z)=h(z)f(z)$ is a meromorphic
function without poles. By Theorem \ref{prod} and Corollary \ref{anal},
we obtain that $g(z)$ is an analytic a.p.\,function.
So we have $f=g/h$.
\hfill\rule{7pt}{7pt}

\begin{cor}
\label{sum}
Let $f_1(z),\ f_2(z)$ be meromorphic a.p.\,functions.
A necessary and sufficient conditions for the sum $f_1(z)+f_2(z)$ to be
a meromorphic a.p.\,function is that zeros and poles of this sum
be separated inside $S$.
\end{cor}

\D\ It follows from Theorem \ref{quot_2} that the sum $f_1(z)+f_2(z)$
is a quotient of two analytic a.p.\,functions; so the assertion
follows from Corollary \ref{quot_1}.
\hfill\rule{7pt}{7pt}

\medskip
Consider the problem of realizability of an a.p.\,divisor as the divisor
of a meromorphic a.p.\,function. This problem was solved in \cite{F1,F2}
for positive a.p.\,divisors and analytic a.p.\,functions by the methods
of cohomology theory. Namely, it was proved that to each positive
a.p.\,divisor $d$ in a strip $S$ a class of \v Cech cohomology
$c(d)$ of the group $H^2(K_{\R},\Z)$ is assigned, $K_{\R}$ being the
universal Bohr
compactification of $\R$;\  $c(d)=0$ iff $d$ is the divisor of an
a.p.\,function on $S$, and $c(d_1+d_2)=c(d_1)+c(d_2)$. Moreover, the
element $c(d)$ remains the same for the restriction of $d$ to any
$S'\subset S$ and for the image of $d$ under every homeomorphism
of $S$ onto $\tilde S$ of the form
$$
\Gamma(x+iy)=x+i\gamma(y).                         \eqno(4)
$$

We do not give the definition of the Bohr
compactification here (one can find it in \cite{D}).
We need only that the group $H^2(K_{\R},\Z)$ can be identified
with the factor group $\R\wedge_{\Z}\R=\R\otimes_{\Z}\R/
\{a\otimes a:\,a\in \R\}$ (see \cite{HM}).
For example, the element $c(d)=\lambda\wedge\mu\in\R\wedge_{\Z}\R$
corresponds to the a.p.\,divisor $d^{\lambda\mu}$ with the support
$\{(\lambda+i\mu)^{-1}n_1+i(\lambda+i\mu)^{-1}n_2:\,n_1,n_2\in \Z\}$
with the multiplicities $k_{n_1,n_2}=1$. If $\lambda/\mu$ is rational,
then $c(d)=0$; otherwise, $c(d^{\lambda\mu })\ne 0$ and
the divisor $d^{\lambda\mu}$ is the divisor of no analytic a.p.\,function;
this coincides with the corresponding result of \cite{T2}.

Note that every element of $\R\wedge_{\Z}\R$, in particular
$c(d)$, is a finite sum of elements $\lambda\wedge\mu$ with irrational
$\lambda/\mu$.

\begin{th}
\label{div_mer}
A divisor $d$ in a strip $S$ is the divisor of a meromorphic a.p.\,function
on $S$ if and only if the following conditions are fulfilled:\\
i) the supports of the divisors $d^{+}$ and $d^{-}$ are separated
inside $S$,\\
ii) $d^+,d^-$ are a.p.\,divisors,\\
iii) $c(d^+)=c(d^-).$
\end{th}

\D\ If $d=d_f$ for some meromorphic a.p.\,function
on $S$, then i) follows from Corollary \ref{sep}, ii) follows from
Theorem \ref{div}. Now by Theorem \ref{quot_2}  we have $f=g/h$,
where $g,\,h$ are analytic a.p.\,function on $S$, therefore
$c(d_g)=c(d_h)=0$. Since $d_{g}= d_{f}^{+}+d',\  d_{h}= d_{f}^{-}+d'$
for any a.p.\,divisor $d'$ in $S$, we obtain iii).

On the other hand, suppose i), ii), iii) are fulfilled; then $c(d^{+})=
\sum\limits_{j=1}^n \lambda_j\wedge\mu_j$ with irrational
$\lambda_j/\mu_j,\ 1\le j\le n$. By $d'$ denote the sum
$\sum_{j=1}^n d^{\mu_j\lambda_j} $; then $ c(d^{+}+d')=
c(d^{+})+\sum\limits_{j=1}^n \mu_j\wedge\lambda_j=0$, hence there exists
an analytic a.p.\,function $g(z)$ on $S$ such that $d_{g}= d_{f}^{+}+d'$;
By iii) we have $c(d^{-}+d')=c(d^++d')=0$, therefore $d^{-}+d'=d_h$ for
some analytic a.p.\,function $h(z)$ on $S$. It follows from i) that
$f=g/h$ is a meromorphic a.p.\,function on $S$.
Finally, $d_f=d_g-d_h=d$.
\hfill\rule{7pt}{7pt}

\begin{cor}
\label{prop}
Suppose $d$ is a divisor in $S$ such that conditions i) and ii) of
Theorem \ref{div_mer} are fulfilled.
If the restriction of $d$ to $S'\subset S$ is the divisor of some
meromorphic a.p.\,function on $S'$, then $d$ is the divisor of a
meromorphic a.p.\,function on $S$; if $\Gamma$ is a homeomorphism of $S$ onto
$\tilde S$ of the form (4) and $d$ is the divisor of some
meromorphic a.p.\,function on $S$, then $\Gamma d$ is the divisor of some
meromorphic a.p.\,function on $\tilde S$.\hfill\rule{7pt}{7pt}
\end{cor}

\bigskip

Department of Mathematics

Kharkov State University

4, pl. Svobody, Kharkov 61077

Ukraine

\newpage

 N.\,Parfyonova, S.\,Favorov.
{\bf Meromorphic almost periodic functions.}
 We introduce a notion of  meromorphic almost periodic function
 and study properties of this class of functions. In particular, we find a
 criterion
 for the product of meromorphic almost periodic functions to be a meromorphic
 almost
 periodic function, too. We prove that every meromorphic almost
 periodic function is a quotient of two analytic almost periodic
 functions.

\bigskip

Н.Д.\,Парф\"енова, С.Ю.Фаворов.
{\bf Мероморфные почти периодичес\-кие функции.}
Введено понятие мероморфной почти периодической функции, изучены свойства
функций этого класса. В частности, найдено необходимое и достаточное
условие того, что произведение двух мероморфных почти периодических
функций опять является мероморфной почти периодической функцией.
Доказано, что любая мероморфная почти периодическая функция является
отношением двух аналитических почти периодических функций.

\bigskip

Н.Д.\,Парфьонова, С.Ю.Фаворов.
{\bf Мероморфнi майже перiодичнi функцii.}
Введено поняття мероморфноi майже перiодичноi функцii, вивченi властивостi
функцiй цього класу. Зокрема, знайдено необхiдна i достатня умова того, що
добуток двух мероморфних майже перiодичних функцiй знов э мероморфною майже
перiодичною функцiэю. Доведено, що будь-яка мероморфна майже перiодична
функцiя э вiдношенням двох аналiтичних майже перiодичних функцiй.

\end{document}